\documentclass[reqno]{amsproc}
\usepackage{graphic x}
\usepackage{amssymb}
\usepackage[margin=1in]{geometry}
\usepackage{setspace,fullpage}
\usepackage{epstopdf}
\usepackage{amsmath,amsthm,amscd,amssymb}
\usepackage{latexsym}
\usepackage[colorlinks,citecolor=red,pagebackref,hypertexnames=false]{hyperref}
\usepackage{geometry}
\numberwithin{equation}{section}
\theoremstyle{plain}
\newtheorem{theorem}{{\bf Theorem}}[section]

\newtheorem{corollary}[theorem]{Corollary}
\newtheorem{proposition}[theorem]{Proposition}
\newtheorem{conjecture}[theorem]{Conjecture}

\theoremstyle{definition}

\newtheorem{claim}[theorem]{Claim}

\newtheorem{case[theorem]}{Case}

\theoremstyle{remark}
\newtheorem{remark}[theorem]{Remark}
\numberwithin{equation}{section}

\begin{document}

\title{\parbox{14cm}{\centering{A note on the multiplicative structure of an additively shifted product set, $AA+1$.}}}

\author{Steven Senger}

\subitem \email{senger@math.udel.edu}

\thanks{}

\maketitle

\setstretch{1.25}
\begin{abstract} We consider the multiplicative structure of sets of the form $AA+1$, where where $A$ is a large, finite set of real numbers. In particular, we show that the additively shifted product set, $AA+1$ must have a large part outside of any generalized geometric progression of comparable length. We prove an analogous result in finite fields as well.
\end{abstract}

\section{Introduction}
There are many problems in additive combinatorics which seek to differentiate between additive and multiplicative structure. By additive (resp. multiplicative) structure in a set, we refer to some arrangement or information that is largely undisturbed by addition (resp. multiplication). A prime example is the sums and products problem. Let $A$ be a large, finite set of $n$ natural numbers. Define the sum set of $A$ to be 
$$A+A = \{ a+b : a\in A, b\in A\}.$$
Define the product set of $A$ to be
$$AA = \{ ab : a\in A, b\in A\}.$$
Let $|\cdot |$ denote the size of a set. The sums and products problem conjectures that,
$$\max\{|A+A|,|AA|\} \geq cn^{2-\epsilon},$$
for any $\epsilon > 0,$ and some constant $c$ which is independent of $n$. In \cite{Elek}, Elekes made progress with an elegant proof based on the celebrated Szemer\'edi-Trotter point-line incidence theorem, from \cite{ST83}. Many of the results in this area have seen largely geometric proofs, including the current record in \cite{Soly08}. Therein, Solymosi proved the remarkable result that either the set of sums or the set of products must have more than about $n^\frac{4}{3}$ elements. See the book by Nathanson, \cite{Nath96}, or the book by Tao and Vu, \cite{TV06}, for more on these and related areas of exploration.

One indication of multiplicative structure is the how the size of a finite set $A$ compares to the size of $AA$. If $AA$ is not much bigger that $A$, then there must be some multiplicative structure in the set $A$. However, regardless of the multiplicative structure in the initial set, the product set will have more. This can be estimated using tools such as the Pl\" unnecke-Rusza inequalities. Again, both \cite{Nath96} and \cite{TV06} provide a good treatment of these results. Our main focus is to show that the multiplicative structure inherent in the product set of a large, finite set of numbers cannot be maximal after an additive shift. This will be made precise in the statement of the main theorem. First, however, we need to introduce some definitions and notation.

In what follows, we use the following asymptotic notation. If two quantities, $X$ and $Y$, vary with respect to some parameter, $n$, we say $X\lesssim Y$ if $X\le CY$, for some constant, $C>0$, which does not depend on $n$. We write $X \approx Y$ when $X\lesssim Y$ and $Y\lesssim X$. If $A$ is a set of numbers, then define its \textit{additive shift} to be
$$A+1 = \{a+1:a\in A\}.$$
Similarly, a \textit{scaling} of $A$ by some number $s$ will be
$$sA=\{sa:a\in A\}.$$
We note that multiplicative behavior of additive shifts have been studied in relation to product sets by Garaev and Shen in \cite{GS}, and Jones and Roche-Newton in \cite{JR}. However, they consider sets of the form $A(A+1)$, which exhibit behavior which is quite different from that of sets of the form $AA+1$, considered here.

Let $r_0, r_1, \dots, r_d$ be real numbers called \textit{generators}, and let $l_1, \dots, l_d$ be positive integers greater than 2. We define the \textit{$d$-dimensional generalized arithmetic progression}
$$R = R(r_0,r_1, \dots, r_d,l_1,\dots, l_d) = \left\{ r_0+x_1 r_1+\dots +x_d r_d : x_j\in \mathbb{Z}, 0\leq x_j <l_j, j=1, \dots, d\right\}.$$
The related notion of a \textit{$d$-dimensional generalized geometric progression} is defined as
$$G=G(g_0,R) = \{g_0^r:r\in R\},$$
where $g_0$ is some positive real number and $R$ is some $d$-dimensional arithmetic progression. We will call either type of generalized progression \textit{degenerate} if $d$ grows asymptotically with the size of the progression. That is to say, if the number of generators is not like a constant compared to the length of the progression, it is degenerate.

Such progressions exhibit maximality in arithmetic, or, respectively, geometric structure. We clarify this with the following elementary proposition.
\begin{proposition}\label{gp} If $R$ is a non-degenerate generalized arithmetic progression, we have that
$$|R+R|\approx |R|.$$
Also, if $G$ is a non-degenerate generalized geometric progression, we have that
$$|GG|\approx |G|.$$
\end{proposition}


We now state the main result.

\begin{theorem}\label{main}
Let $A \subset \mathbb{R}$, be a large, finite set of numbers. Let $G$ be any non-degenerate generalized geometric progression with $|G| \approx |AA|$. We have that
$$|(AA+1) \setminus G| \gtrsim |A|^{1-\delta}$$
for any $\delta >0$.
\end{theorem}

We remark that the proof of Theorem \ref{main} holds true for some slightly more general progressions which we consider degenerate. See the remark following Claim \ref{bb}. Two direct corollaries follow.

\begin{corollary}
Let $A \subset \mathbb{R}$, be a large, finite set of numbers. Let $G$ be any non-degenerate generalized geometric progression with $|G| \approx |AA|$. We have that
$$(AA+1) \not\subset G.$$
\end{corollary}
Notice that for any non-degenerate geometric progression, $H\subset \mathbb{R}$, there exists a set $H'\subset \mathbb{R}$ such that $|H'H'|\approx |H|$ and $H \subset H'H'.$ If we apply Theorem \ref{main} with $A=H'$, we get the following corollary.
\begin{corollary}
Let $G$ and $H$ be any two large, finite, non-degenerate generalized geometric progressions with $|G| \approx |H|$. We have that
$$(H+1) \not\subset G.$$
\end{corollary}

We suspect that the additive shift disrupts multiplicative structure even more than Theorem \ref{main} indicates, as suggested by the following conjecture.

\begin{conjecture}
Let $A \subset \mathbb{R}$, be a large, finite set of numbers. If
$$|(AA+1)\cap BC| \approx |AA|,$$
where $B$ and $C$ are also large finite sets of numbers, then $\min\{|B|,|C|\}\lesssim 1$.
\end{conjecture}

The next result is of a similar type, but in the setting of finite fields.

\begin{theorem}\label{ffmain}
Let $A\subset \mathbb{F}_q$ such that the following two conditions hold:
\begin{enumerate}
\item There exists a real number $\epsilon > 0$ such that $|A||AA|\gtrsim q^{\frac{3}{2}+ \epsilon}.$
\item There exists a real number $\delta >0$ such that $|AA|\lesssim q^{1-\delta}.$
\end{enumerate}
Let $G$ be any non-degenerate generalized geometric progression with $|G| \approx |AA|$. We have that
$$|(AA+1) \setminus G| \gtrsim q^\delta.$$
\end{theorem}


\section{Proof of Theorem \ref{main}}
The basic outline of the proof is to start with a given large, finite subset of $\mathbb{R}$. Then, with this set and any appropriate generalized geometric progression, we construct two large, finite sets of points in $\mathbb{R}^2$. We will then apply the a recent result by Iosevich, Roche-Newton, and Rudnev regarding the set of dot products determined by our point sets, from \cite{IRR}. The underlying arithmetic of the dot product set will allow us to compare the elements of the shifted product set to the elements of the progression.

The key ingredient to their proof is inspired by recent developments in the study of the Erd\H os distance problem. The classical Erd\H os distance problem asks for the minimum number of distinct distances which can be determined by any large, finite set of $n$ points. The conjecture in the plane was that any such set must determine at least $n^{1-\epsilon}$ distinct distances, for any $\epsilon>0.$ Guth and Katz proved this in \cite{GK10} with a blend of cell-decomposition and algebraic geometry, applied to an incidence problem in three dimensions. Shortly thereafter, Iosevich, Roche-Newton, and Rudnev, used similar techniques to prove a related result on the number of distinct dot products determined by such point sets in the plane. Specifically, they proved the following theorem.

\begin{theorem} \label{dp}
Consider any large finite point sets $E,F\subset \mathbb{R}^2$ of $n$ points each, neither of which is contained in a single line. Let $\Pi(E,F)$ denote the set of dot products
$$\Pi(E,F) = \{x\cdot y:x\in E, y\in F\}.$$
Then, for any $\epsilon>0,$ the number of distinct dot products is bounded below by
$$|\Pi(E,F)| \gtrsim n^{1-\epsilon}.$$
\end{theorem}

We now prove Theorem \ref{main}.

\begin{proof}
Fix any large finite set $A\subset \mathbb{R}$, and a real number $\delta >0$. Let $G=G(g_0,R)$ be any non-degenerate generalized geometric progression with $|G|\approx |AA|$. Consider $R=R(r_0,r_1, \dots, r_d,l_1,\dots, l_d)$, the $d$-dimensional arithmetic progression defining the exponents of $g_0$ which make up $G$. Since $G$ is non-degenerate, $R$ must also be non-degenerate. We will define $g_1$ to be the ``first element" of $G$, namely,
$$g_1:=g_0^{r_0}.$$
In what follows, we need to work with the normalized progression, $G'$ which will be defined as
$$G':=\frac{G}{g_1} = \left\lbrace \frac{g}{g_1}: g\in G\right\rbrace.$$
Notice that $|G|=|G'|$. Now, define the set $B$ to be
$$B:=\lbrace g\in G': gg\in G'\rbrace = \lbrace g\in G': g=g_0^{x_1 r_1+\dots +x_d r_d}, x_j < l_j,x_j\in 2\mathbb{N}, j=1, \dots, d \rbrace.$$

\begin{claim}\label{bb}
$|B|\approx|G|$.
\end{claim}
Since $B \subset G'$, it is clear that $|B| \lesssim |G|$. Now we need only to show that $|B|\gtrsim |G|$. Notice that for an element to be in $B$, its square must be in $G'$, hence the evenness condition on the $x_j$ in the definition of $B$. So, we can count the number of elements in $B$ by counting the number of elements of $G'$ whose corresponding $x_j$ are all even. By definition of the $l_j$, we get
$$|B|\geq \prod_{j=1}^d \left\lfloor \frac{l_j}{2} \right\rfloor \geq \prod_{j=1}^d \frac{l_j}{3}\geq \frac{|G|}{3^d} ,$$
as $|G|$ is equal to the product of the $l_j$, and the claim is proved.


By Claim \ref{bb}, $|B|\approx |G|$, so we have that $|BB|\approx |AA|.$
Now, we construct $E$ and $F$, finite subsets of $\mathbb{R}^2$,
$$E:=\left\{\left(bg_1,ba\right)\in \mathbb{R}^2 :b\in B, a\in A\right\}, \text{and}~F:=\left\{\left(b,ba\right)\in \mathbb{R}^2 :b\in B, a\in A\right\}.$$
These sets will have size $|E|=|F|=|A||B|\approx|A||G|.$

\begin{remark}
Note that the actual size estimate in Claim \ref{bb} could be satisfied by some slightly more general types of progressions. That is, as long as the dependence of $d$ on the length of the progression is sub-logarithmic, the proof still works.
\end{remark}

We will consider $\Pi(E,F)$, the set of distinct dot products determined by pairs in $E\times F$. Notice that
\begin{align*}
\Pi(E,F) &= \left\{\left(g_1b,ba\right)\cdot\left(b',b'a'\right):a,a'\in A, b,b' \in B\right\}\\
&= \left\{g_1bb'(aa'+1):a,a'\in A, b,b' \in B\right\}\\
&= g_1BB(AA+1)\subset G(AA+1).
\end{align*}
By construction,
$$|\Pi(E,F)|=|g_1BB(AA+1)| \leq |G(AA+1)|.$$

Set $\epsilon = \delta/3$. Since $\epsilon>0$, Theorem \ref{dp} gives us that
$$ |\Pi(E,F)| \gtrsim |E|^{1-\epsilon}\gtrsim (|A||B|)^{1-\epsilon}.$$
Comparing upper and lower bounds on $|\Pi(E,F)|$ gives us
\begin{equation} \label{ABep}
|G(AA+1)|\gtrsim (|A||B|)^{1-\epsilon}.
\end{equation}
Our aim is to get a lower bound on the exceptional set $C := (AA+1) \setminus G$. From \eqref{ABep}, we get
\begin{align*}
(|A||B|)^{1-\epsilon} &\lesssim |G(AA+1)|\\
&=|G((G\cap(AA+1))\cup C)|\\
&=|G(G\cap(AA+1))\cup GC|\\
\end{align*}
Notice that the first term in the above union is a subset of $GG$, and therefore has size $\lesssim |G|$, by Proposition \ref{gp}. So we can conclude that
$$|GC| \gtrsim (|A||B|)^{1-\epsilon}\approx (|A||G|)^{1-\epsilon},$$
which, by simple counting, gives us that
\begin{equation}\label{cag}
|C|\gtrsim \frac{(|A||G|)^{1-\epsilon}}{|G|}= |A|^{1-\epsilon}|G|^{-\epsilon}.
\end{equation}
Again, by a simple counting argument, we see that $|A|^2\gtrsim|AA|$. Since $|G|\approx |AA|$, we can rewrite \eqref{cag} as
\begin{align*}
|C| &\gtrsim |A|^{1-\epsilon}|G|^{-\epsilon} \\
&\gtrsim |A|^{1-\epsilon}|A|^{-2\epsilon} \\
&\gtrsim |A|^{1-\delta},
\end{align*}
where the last line follows by definition of $\epsilon$.

\end{proof}

\section{Proof of Theorem \ref{ffmain}}
We follow a similar program to the proof of Theorem \ref{main}. Therefore, some details are omitted. The dot product set estimate in the finite field setting is in a slightly different form. It is due to Hart, Iosevich, Koh, and Rudnev. The statement of the theorem in \cite{HIKR} is for one set, but the proof works, with obvious modifications, for two sets. We consider the special case that both sets have the same size. For this section, define $\Pi(E,F)$ as before, except for subsets of $\mathbb{F}_q^2$ instead of $\mathbb{R}^2$. Also, let $\mathbb{F}_q^*$ denote the multiplicative group of $\mathbb{F}_q$.

\begin{theorem}\label{ffdp}
Let $E,F \subset \mathbb{F}_q^d$ such that $|E|=|F|>q^\frac{d+1}{2}$. Then
$$\mathbb{F}_q^*\subset\Pi(E,F) = \{x\cdot y:x\in E, y\in F\}.$$
\end{theorem}

Notice that the size condition in Theorem \ref{ffdp} is given with constant 1. This is why we include the $\epsilon$ in the size condition of Theorem \ref{ffmain}, although a slightly more general statement is also true. Also, notice that we only use the case that $d=2$.

\begin{proof}
Let $A\subset\mathbb{F}_q$ be given, and suppose that it satisfies the two size conditions:
\begin{enumerate}
\item There exists a real number $\epsilon > 0$ such that $|A||AA|\gtrsim q^{\frac{3}{2}+\epsilon}.$
\item There exists a real number $\delta >0$ such that $|AA|\lesssim q^{1-\delta}.$
\end{enumerate}
Now, let $G=G(g_0,R)$ be any non-degenerate generalized geometric progression with $|G|\approx |AA|.$ Again, we define
$$g_1:=g_0^{r_0}, G':=\frac{G}{g_1} = \left\lbrace \frac{g}{g_1}: g\in G\right\rbrace, \text{ and } B:=\lbrace g\in G': gg\in G'\rbrace.$$
Claim \ref{bb} still holds in this context, so, as before, we have $|B|\approx |BB|\approx |AA|.$
Now, we construct $E$ and $F$, finite subsets of $\mathbb{F}_q^2$,
$$E:=\left\{\left(bg_1,ba\right)\in \mathbb{F}_q^2 :b\in B, a\in A\right\}, \text{and}~F:=\left\{\left(b,ba\right)\in \mathbb{F}_q^2 :b\in B, a\in A\right\}.$$
As in the proof of Theorem \ref{main}, the set of dot products determined by pairs in $E\times F$ will be
$$\Pi(E,F) = g_1BB(AA+1).$$
We also know that $|E|=|F|=|A||B|\approx|A||G|.$ So, by the first size condition satisfied by $A$, and the fact that $|G| \approx |AA|$, we see that $|E|\gtrsim q^{\frac{3}{2}+\epsilon}$. Since $E$ is large enough to satisfy the hypotheses of Theorem \ref{ffdp}, we are guaranteed that $|\Pi(E,F)|\geq q-1$. Specifically, using the proof of Theorem \ref{main} as a guide, we get
\begin{align}
\nonumber q-1 &\leq |\Pi(E,F)|\\ 
\nonumber &=|g_1BB(AA+1)|\\ 
\label{q} &\leq |G(AA+1)|. 
\end{align}
We again seek a lower bound on the exceptional set. Define $C \subset \mathbb{F}_q$ to be $(AA+1) \setminus G$. By \eqref{q} and the definition of $C$, we get
\begin{align*}
q-1 &\leq |G(AA+1)|\\
&=|G((G\cap(AA+1))\cup C)|\\
&=|G(G\cap(AA+1))\cup GC|\\
\end{align*}
Again, the first term in the union will have size $\lesssim |G|\approx |AA|$. The second size condition satisfied by $A$ tells us that $|AA|\lesssim q^{1-\delta}$, so the second term dominates. This gives us that $|GC| \approx q$, which, by simple counting and the fact that $|G|\approx |AA|$, yields
$$|C|\gtrsim \frac{q}{|AA|}=q^\delta,$$
as claimed.
\end{proof}


\vskip.5in

\end{document}